\documentclass[10pt,draft]{amsart}
\usepackage{amscd}
\usepackage{amssymb}

\newtheorem{thm}[equation]{Theorem}
\newtheorem{cor}[equation]{Corollary}
\newtheorem{prop}[equation]{Proposition}
\newtheorem{lem}[equation]{Lemma}
\newtheorem{defn}[equation]{Definition}
\newtheorem{rem}[equation]{Remark}
\newtheorem{exmp}[equation]{Example}

\numberwithin{equation}{section}

\newcommand{\arr}{\rightarrow}
\newcommand{\larr}{\leftarrow}
\newcommand{\lgarr}{\longrightarrow}

\newcommand{\xarr}{\xrightarrow}

\newcommand{\cat}[1]{\operatorname{\mathsf{#1}}}

\newcommand{\opn}{\operatorname}

\newcommand{\rmitem}[1]{\item[\text{\textup{(#1)}}]}

\newcommand{\mcal}[1]{\mathcal{#1}}

\newcommand{\mrm}[1]{\mathrm{#1}}
\newcommand{\mbb}[1]{\mathbb{#1}}

\title[On $t$-structures and Torsion Theories]{On $t$-structures and 
Torsion Theories \\ Induced by Compact Objects}
\author{Mitsuo Hoshino, Yoshiaki Kato and Jun-ichi Miyachi}
\date{1 May 2000}
\address{M. Hoshino: Department of Mathematics, University of Tsukuba,
Ibaraki, 305-8571, Japan}
\email{hoshino@math.tsukuba.ac.jp}
\address{Y. Kato: Department of Mathematics, University of Tsukuba,
Ibaraki, 305-8571, Japan}
\email{ykato@math.tsukuba.ac.jp}
\address{J. Miyachi: Department of Mathematics, Tokyo Gakugei
University, Koganei-shi, Tokyo, 184-8501, Japan}
\email{miyachi@u-gakugei.ac.jp}

\begin{document}

\begin{abstract}
First, we show that a compact object $C$ in a triangulated category,
which satisfies suitable conditions, induces a $t$-structure.
Second, in an abelian category we show that a complex $P^{\centerdot}$ 
of small projective objects of term length two,
which satisfies suitable conditions, induces a torsion theory. 
In the case of module categories, using a torsion theory, we give 
equivalent conditions for $P^{\centerdot}$ to be a tilting complex.
Finally, in the case of artin algebras, we give
one to one correspondence between tilting complexes of term length two
and torsion theories with certain conditions.
\end{abstract}

\maketitle

\setcounter{section}{-1}
\section{Introduction} \label{s0}

In the representation theory of finite dimensional algebras, torsion theories were 
studied by several authors in connection with classical tilting modules.  
For these torsion theories, there are equivalences between torsion (resp., torsionfree) 
classes and torsionfree (resp., torsion) classes, which is known as 
Theorem of Brenner and Butler (\cite{HR}).
One of the authors gave one to one correspondence between 
classical tilting modules and torsion theories with certain conditions (\cite{Ho1},
\cite{Ho2}).
But in the case of a self-injective algebra $A$, tilting modules are essentially 
isomorphic to $A$.
In \cite{Ri}, Rickard introduced the notion of tilting complexes as a generalization of 
tilting modules, 
and showed that these complexes induce equivalences between derived categories
of module categories.  Tilting complexes of term length two are often studied
in the case of self-injective algebras (e.g. \cite{Hl}, \cite{HK}).
On the other hand, for triangulated categories, Beilinson, Bernstein and Deligne introduced 
the notions of $t$-structures and admissible abelian subcategories, and studied relationships 
between them (\cite{BBD}).
In this paper, first, we deal with a compact object $C$ in a triangulated category, 
and study a $t$-structure induced by $C$.
Second, in an abelian category $\mcal{A}$ we deal with a complex 
$P^{\centerdot}$ of small projective objects of term length two and 
study a torsion theory induced by $P^{\centerdot}$.

In Section \ref{s1}, we show that a compact object $C$ in a triangulated category $\mcal{T}$, 
which satisfies suitable conditions, induces a $t$-structure
$(\mcal{T}^{\leq 0}(C), \mcal{T}^{\geq 0}(C))$, and its core
$\mcal{T}^{0}(C)$ is equivalent to the category $\cat{Mod}B$
of left $B$-modules,
where  $B = \opn{End}_{\mcal{T}}(C)^{\mrm{op}}$ (Theorem \ref{thm1.3}).
In Section \ref{s2}, we define subcategories $\mcal{X}(P^{\centerdot})$, 
$\mcal{Y}(P^{\centerdot})$ of  an abelian category $\mcal{A}$ satisfying the condition Ab4, and
show when $(\mcal{X}(P^{\centerdot}), \mcal{Y}(P^{\centerdot}))$ is a torsion
theory (Theorem \ref{thm2.10}).  Furthermore, we show that if $P^{\centerdot}$ 
induces a torsion theory
$(\mcal{X}(P^{\centerdot}), \mcal{Y}(P^{\centerdot}))$ for $\mcal{A}$, then the core 
$\cat{D}(\mcal{A})^{0}(P^{\centerdot})$
is admissible abelian, and then there is a torsion theory
$(\mcal{Y}(P^{\centerdot})[1], \mcal{X}(P^{\centerdot}))$
for $\cat{D}(\mcal{A})^{0}(P^{\centerdot})$ (Theorem \ref{thm2.15}).
In Section \ref{s3}, we apply results of  Section \ref{s2} to module categories.
We  characterize a torsion theory for the category $\cat{Mod}A$
of left $A$-modules, and for its core  $\cat{D}(\cat{Mod}A)^{0}(P^{\centerdot})$ 
(Theorems \ref{thm3.5} and \ref{thm3.8}).  Furthermore, using a torsion theory,
we give equivalent conditions for $P^{\centerdot}$ to be a tilting complex 
(Corollary \ref{cor3.6}).
In Section \ref{s4}, We show that, if  $P^{\centerdot}$ is a tilting complex, then
it induces  equivalences between torsion theories for $\cat{Mod}A$ and
for $\cat{Mod}B$, where 
$B = \opn{End}_{\cat{D}(\cat{Mod}A)}(P^{\centerdot})^{\mrm{op}}$ (Theorem \ref{thm4.4}).
In Section \ref{s5}, in the case of artin algebras, if  a torsion theory 
$(\mcal{X}, \mcal{Y})$ satisfies certain conditions, then there exists a
tilting complex $P^{\centerdot}$ of term length two such that a torsion theory  $(\mcal{X}, \mcal{Y})$ 
coincides with $(\mcal{X}(P^{\centerdot}), \mcal{Y}(P^{\centerdot}))$ 
(Theorem \ref{thm5.8}).  As a consequence,
we have one to one correspondence between tilting complexes of term
length two and torsion theories with certain conditions 
(Corollary \ref{cor3.7}, Propositions \ref{prop5.5}, \ref{prop5.7} and Theorem \ref{thm5.8}).

\section{$t$-structures Induced by Compact Objects} \label{s1}

In this section, we deal with a triangulated category $\mcal{T}$ and its 
full subcategory $\mcal{C}$.  We will call $\mcal{C}$ admissible abelian 
provided that $\opn{Hom}_{\mcal{T}}(\mcal{C}, \mcal{C}[n]) = 0$
for $n < 0$, and that all morphisms in $\mcal{C}$ are $\mcal{C}$-admissible 
in the sense of \cite{BBD}, 1.2.3.  In this case, according to \cite{BBD}, 
Proposition 1.2.4, $\mcal{C}$ is an abelian category.  A triangulated 
category $\mcal{T}$ is said to contain direct sums 
if direct sums of objects indexed by  any set exist in $\mcal{T}$.
An object $C$ of $\mcal{T}$ is called compact if 
$\opn{Hom}_{\mcal{T}}(C, -)$ commutes with direct sums.  
Furthermore, a collection $\mcal{S}$ of compact objects of $\mcal{T}$ is called 
a generating set provided that $X = 0$ whenever 
$\opn{Hom}_{\mcal{T}}(\mcal{S}, X) = 0$, and that $\mcal{S}$ is stable 
under suspension (see \cite{Ne} for details).  For an object $C \in \mcal{T}$
and an integer $n$, we denote by 
$\mcal{T}^{\geq n}(C)$ (resp., $\mcal{T}^{\leq n}(C)$)
the full subcategory of $\mcal{T}$ consisting of $X \in \mcal{T}$
with $\opn{Hom}_{\mcal{T}}(C, X[i]) = 0$ for $i < n$ (resp., $i > n$),
and set $\mcal{T}^{0}(C) = \mcal{T}^{\leq 0}(C) \cap \mcal{T}^{\geq 0}(C)$.

For an abelian category $\mcal{A}$, we denote by 
$\cat{C}(\mcal{A})$ the category of complexes of $\mcal{A}$, and denote by 
$\cat{D}(\mcal{A})$ (resp., $\cat{D}^{+}(\mcal{A})$, $\cat{D}^{-}(\mcal{A})$, 
$\cat{D}^{\mrm{b}}(\mcal{A})$)
the derived category of complexes of $\mcal{A}$ (resp.,
complexes of $\mcal{A}$ with bounded below homologies, 
complexes of $\mcal{A}$ with bounded above homologies, 
complexes of $\mcal{A}$ with bounded homologies).
For an additive category $\mcal{B}$, we denote by 
$\cat{K}(\mcal{B})$ (resp., $\cat{K}^{-}(\mcal{B})$, $\cat{K}^{\mrm{b}}(\mcal{B})$)
the homotopy category of complexes of $\mcal{B}$ (resp.,
bounded above complexes of $\mcal{B}$, 
bounded complexes of $\mcal{B}$) (see \cite {RD} for details). \par\noindent

\begin{prop} \label{prop1.1}
Let $\mcal{T}$ be a triangulated category which 
contains direct sums, $C$ a compact object satisfying 
$\opn{Hom}_{\mcal{T}}(C, C[n]) = 0$ for $n > 0$.
Then for any $r \in \mbb{Z}$ and any object $X \in \mcal{T}$,
there exist an object  $X^{\geq r} \in \mcal{T}^{\geq r}(C)$ and a morphism 
${\alpha}^{\geq r} : X \arr X^{\geq r}$ in $\mcal{T}$ such that
\begin{enumerate}
\rmitem{i} for any $i \geq r$,
$\opn{Hom}_{\mcal{T}}(C, {\alpha}^{\geq r} [i])$ is an isomorphism,
\rmitem{ii} for every object $Y \in \mcal{T}^{\geq r}(C)$, 
$\opn{Hom}_{\mcal{T}}({\alpha}^{\geq r} , Y)$ is an isomorphism.
\end{enumerate}
\end{prop}

\begin{proof}  Let $X_{0} = X$.  For $n \geq 1$, by induction we construct 
a distinguished triangle 
\[
C_{n}[n-r] \xarr{g_n} X_{n-1} \xarr{h_n} X_{n} \arr 
\]
as follows. If $\opn{Hom}_{\mcal{T}}(C, X_{n-1}[r-n]) = 0$,  then we set $C_{n} = 0$.  
Otherwise, we take a direct sum $C_{n}$ of copies of $C$ and a morphism 
$g'_n : C_n  \arr X_{n-1}[r-n]$ such 
that $\opn{Hom}_{\mcal{T}}(C,  g'_n)$ is an epimorphism, and let 
$g_n = g'_n[n-r]$.
Then, by easy calculation, we have the following:
\[ \begin{aligned}
(a)\qquad
& \opn{Hom}_{\mcal{T}}(C,  X_n[i]) = 0 \enskip\text{for}\enskip r-n \leq i < r , \\
(b)\qquad
& \opn{Hom}_{\mcal{T}}(C, h_n[i]) \enskip\text{is an isomorphism for any} 
\enskip n \enskip\text{and}\enskip i \geq r.
\end{aligned}\]
Let $X^{\geq r}$ be a homotopy colimit $\underset{\lgarr}{\cat{hocolim}}\enskip X_n$
and ${\alpha}^{\geq r} : X \arr X^{\geq r}$ a structural morphism $X_0 \arr
\underset{\lgarr}{\cat{hocolim}}\enskip X_n$.  According to \cite{Ne}, Lemma 2.8, 
the conditions 
($a$), ($b$) imply that $X^{\geq r}$ belongs to $\mcal{T}^{\geq r}(C)$ and 
satisfies the statement (i).
For an object $Y \in \mcal{T}^{\geq r}(C)$, we have an exact sequence
\[ \begin{aligned}
\opn{Hom}_{\mcal{T}}(C_n[n-r], Y[j-1]) \arr
\opn{Hom}_{\mcal{T}}(X_n, Y[j])  & \arr \\
\opn{Hom}_{\mcal{T}}(X_{n-1}, Y[j]) & \arr
\opn{Hom}_{\mcal{T}}(C_n[n-r], Y[j]) .
\end{aligned} \]
Since $\opn{Hom}_{\mcal{T}}(C[i], Y[j]) = 0$ for $j-i < r$,
$\opn{Hom}_{\mcal{T}}(h_n, Y[j])$ is an isomorphism for any $n \geq 1$ 
and $j \leq 0$.  Then, we have an epimorphism
\[
{\prod}_{n}\opn{Hom}_{\mcal{T}}(X_n, Y[j])
\xarr{1- \text{shift}}
{\prod}_{n}\opn{Hom}_{\mcal{T}}(X_n, Y[j])
\]
for any $j \leq 0$.  Therefore, we have an exact sequence
\[
0 \arr
\opn{Hom}_{\mcal{T}}(X^{\geq r}, Y) \arr
{\prod}_{n}\opn{Hom}_{\mcal{T}}(X_n, Y)
\xarr{1- \text{shift}}
{\prod}_{n}\opn{Hom}_{\mcal{T}}(X_n, Y)
\arr 0 .
\]
Hence we have
\[ \begin{aligned}
\opn{Hom}_{\mcal{T}}(X^{\geq r}, Y) & \cong 
\varprojlim\opn{Hom}_{\mcal{T}}(X_n, Y) \\
& \cong \opn{Hom}_{\mcal{T}}(X, Y) .
\end{aligned} \]
\end{proof}

\begin{defn}[\cite{BBD}] \label{dfn1.2} 
Let $\mcal{T}$ be a triangulated category.  For full subcategories
$\mcal{T}^{\leq 0}$ and $\mcal{T}^{\geq 0}$, $(\mcal{T}^{\leq 0}, \mcal{T}^{\geq 0})$
is called a t-structure on $\mcal{T}$ provided that
\begin{enumerate}
\rmitem{i} $\opn{Hom}_{\mcal{T}}(\mcal{T}^{\leq 0},\mcal{T}^{\geq 1}) = 0$;
\rmitem{ii} $\mcal{T}^{\leq 0}  \subset \mcal{T}^{\leq 1}$ and
$\mcal{T}^{\geq 0} \supset \mcal{T}^{\geq 1}$;
\rmitem{iii} for any $X \in \mcal{T}$, there exists a distinguished triangle
\[
X' \arr X \arr X'' \arr
\]
with $X' \in \mcal{T}^{\leq 0}$ and $X'' \in \mcal{T}^{\geq 1}$,
\end{enumerate}
where $\mcal{T}^{\leq n} = \mcal{T}^{\leq 0}[-n]$ and
$\mcal{T}^{\geq n} = \mcal{T}^{\geq 0}[-n]$.

A t-structure $(\mcal{T}^{\leq 0}, \mcal{T}^{\geq 0})$ on $\mcal{T}$
is called non-degenerate if ${\bigcap}_{n \in \mbb{Z}}\mcal{T}^{\leq n} =
{\bigcap}_{n \in \mbb{Z}}\mcal{T}^{\geq n} \\ = \{ 0\}$.
\end{defn}

\begin{thm} \label{thm1.3}  Let $\mcal{T}$ be a triangulated category which 
contains direct sums, $C$ a compact object satisfying 
$\opn{Hom}_{\mcal{T}}(C, C[n]) = 0$ for $n > 0$, and 
$B = \opn{End}_{\mcal{T}}(C)^{\mrm{op}}$.
If $\{C[i] : i \in \mbb{Z}\}$ is a generating set, then the following hold.
\begin{enumerate}
 \rmitem{1} $(\mcal{T}^{\leq 0}(C), \mcal{T}^{\geq 0}(C))$ is a non-degenerate
 t-structure on $\mcal{T}$.
 \rmitem{2} $\mcal{T}^{0}(C)$ is admissible abelian.
 \rmitem{3} The functor
\[
\opn{Hom}_{\mcal{T}}(C, -) : \mcal{T}^{0}(C) \arr \cat{Mod}B
\]
is an equivalence.
\end{enumerate}
\end{thm}

\begin{proof}
(1)  For any object $X \in \mcal{T}^{\leq 0}(C)$, we take an object $X^{\geq 1}
\in \mcal{T}^{\geq 1}(C)$ and a morphism ${\alpha}^{\geq 1}: X \arr X^{\geq 1}$
satisfying the conditions of Proposition \ref{prop1.1}.  Then
for any $Y \in \mcal{T}^{\geq 1}(C)$, by Proposition \ref{prop1.1} (ii), we have
\[
\opn{Hom}_{\mcal{T}}(X^{\geq 1},Y) \cong \opn{Hom}_{\mcal{T}}(X,Y).
\]
By Proposition \ref{prop1.1} (i), $X \in \mcal{T}^{\leq 0}(C)$ implies that
$\opn{Hom}_{\mcal{T}}(C, X^{\geq 1}[i]) = 0$ for all $i \in \mbb{Z}$.
Since $\{C[i] : i \in \mbb{Z}\}$ is a generating set, we have $X^{\geq 1} = 0$,
and hence $\opn{Hom}_{\mcal{T}}(X,Y) = 0$.
It is easy to see that $\mcal{T}^{\leq 0}(C)  \subset \mcal{T}^{\leq 1}(C)$ and
$\mcal{T}^{\geq 0}(C) \supset \mcal{T}^{\geq 1}(C)$.
For any object $Z \in \mcal{T}$, we take an object $Z^{\geq 1}
\in \mcal{T}^{\geq 1}(C)$ and a morphism ${\alpha}^{\geq 1}: Z \arr Z^{\geq 1}$
satisfying the conditions of Proposition \ref{prop1.1}, and embed ${\alpha}^{\geq 1}$
in a distinguished triangle
\[
Z' \arr Z \arr Z^{\geq 1} \arr .
\]
Applying $\opn{Hom}_{\mcal{T}}(C,-)$ to the above triangle, by Proposition 
\ref{prop1.1} (i), we have $Z' \in \mcal{T}^{\leq 0}(C)$.
Since $\{C[i] : i \in \mbb{Z}\}$ is a generating set, it is easy to see that
$(\mcal{T}^{\leq 0}(C), \\ \mcal{T}^{\geq 0}(C))$ is non-degenerate.

(2)  Since $\mcal{T}^{0}(C)$ is the core of the $t$-structure
$(\mcal{T}^{\leq 0}(C), \mcal{T}^{\geq 0}(C))$, 
the assertion follows by \cite{BBD}, Theorem 1.3.6.

(3)  Step 1:  According to \cite{BBD}, Proposition 1.2.2, 
the short exact sequences in $\mcal{T}^{0}(C)$ are just the distinguished triangles
\[
X  \arr Y \arr Z \arr 
\]
with $X, Y$ and $Z$ belonging to $\mcal{T}^{0}(C)$.  
It follows that $\opn{Hom}_{\mcal{T}}(C, - ) :  \mcal{T}^{0}(C) \arr \cat{Mod}B$
is exact.
Let $M \in \cat{Mod}B$ and take a free presentation $P_1 \arr P_0 \arr M
\arr 0$.  We take $C' = C^{\geq 0} 
\in \mcal{T}^{0}(C)$ and $\alpha = {\alpha}^{\geq 0} : C \arr C'$ satisfying the conditions of 
Proposition \ref{prop1.1}.  
Then there exist sets $I, J$ and a 
collection of morphisms $h_{ij} : C' \arr C'$ such that
\[\begin{CD}
P_1 @>>> P_0 \\
@VVV @VVV \\
\opn{Hom}_{\mcal{T}}(C, C')^{(J)} @>{\oplus}_{ij}\opn{Hom}(C,h_{ij})>>
\opn{Hom}_{\mcal{T}}(C, C')^{(I)}
\end{CD} \]
is commutative, where the vertical arrows are isomorphisms.  We take an 
exact sequence in $\mcal{T}^{0}(C)$
\[
{C'}^{(J)} \xarr{{\oplus}_{ij}h_{ij}} {C'}^{(I)} \arr X \arr 0.
\]
Since $C$ is compact, by the exactness of $\opn{Hom}_{\mcal{T}}(C, -)$, we have  
$\opn{Hom}_{\mcal{T}}(C, X) \cong M$.

\vspace{12pt}
Step 2:  We show that $\opn{Hom}_{\mcal{T}}(C, -)$ reflects isomorphisms.  Let
\[
X \xarr{u} Y \arr Z \arr 
\]
be a distinguished triangle in $\mcal{T}$ with $X, Y \in \mcal{T}^{0}(C)$ and 
with $\opn{Hom}_{\mcal{T}}(C, u)$ an isomorphism.  Then, by applying 
$\opn{Hom}_{\mcal{T}}(C, -)$, we get $\opn{Hom}_{\mcal{T}}(C, Z[n]) = 0$ for all 
$n \in \mbb{Z}$, and hence $Z = 0$.  It follows that $u$ is an 
isomorphism.

Next, we show that $\opn{Hom}_{\mcal{T}}(C, -)$ is faithful.  Let $v: X 
\arr Y$ be a morphism in  $\mcal{T}^{0}(C)$ with  $\opn{Hom}_{\mcal{T}}(C, v) = 0$. 
By the exactness of $\opn{Hom}_{\mcal{T}}(C, -)$, $\opn{Hom}_{\mcal{T}}(C, 
\opn{Im}v) \cong \opn{Im}\opn{Hom}_{\mcal{T}}(C, v) = 0$.  Since $\opn{Im}v 
\in \mcal{T}^{0}(C)$, we have $\opn{Hom}_{\mcal{T}}(C, \opn{Im}v[n]) = 0$ for 
all $n \in \mbb{Z}$, and hence $\opn{Im}v = 0$ and $v = 0$. 

Let $\mcal{M}$ be the full subcategory of $\mcal{T}^{0}(C)$ consisting of objects
$X$ such that there exists an exact sequence $C_1  \arr  C_0 \arr X \arr 0$ in 
$\mcal{T}^{0}(C)$, where $C_0, C_1$ are direct sums of copies of $C'$.
Since $\opn{Hom}_{\mcal{T}}(C, -)$ is faithful, by the same technique 
as in (1), it is not hard to see that $\opn{Hom}_{\mcal{T}}(C, - )|_{\mcal{M}}$ 
is full dense, and hence an equivalence.  It remains to show that 
$\mcal{M} = \mcal{T}^{0}(C)$.  
For an object $X  \in \mcal{T}^{0}(C)$, we have a commutative diagram
\[ \begin{CD}
\opn{Hom}_{\mcal{T}}(C, C^{(J)})  @>\opn{Hom}(C, f)>> \opn{Hom}_{\mcal{T}}(C, C^{(I)})
@>\opn{Hom}(C, g)>> \opn{Hom}_{\mcal{T}}(C, X) @>>> 0 \\
@VV \opn{Hom}_{\mcal{T}}(C, {\alpha}_{1}) V
@VV \opn{Hom}_{\mcal{T}}(C, {\alpha}_{0}) V \\
\opn{Hom}_{\mcal{T}}(C, {C'}^{(J)}) @>\opn{Hom}(C, f')>>  \opn{Hom}_{\mcal{T}}(C, {C'}^{(I)})
\end{CD} \]
with the top row being exact and with the vertical arrows being isomorphisms.
And we have a commutative diagram in $\mcal{T}$
\[ \begin{CD}
C^{(J)} @>f>>  C^{(I)} @>g>>  X \\
@VV {\alpha}_{1} V @VV {\alpha}_{0}V \\
{C'}^{(J)} @>f'>>  {C'}^{(I)}
\end{CD} \]
with $gf = 0$.  By Proposition \ref{prop1.1}(ii), there exists $h: {C'}^{(I)} 
\arr X$ such that $g = h{\alpha}_{0}$.  Since $\opn{Hom}_{\mcal{T}}(C, hf') = 0$,
we have $hf' = 0$.  Then there exists $w : \opn{Cok}f' \arr X$ such that 
$g = w g'{\alpha}_{0}$, where $g': {C'}^{(I)} \arr \opn{Cok}f'$ is a canonical 
morphism.  Then $\opn{Hom}_{\mcal{T}}(C, w )$ is an isomorphism, 
and hence $w$ is an isomorphism and $X \cong \opn{Cok}f' \in \mcal{M}$.
\end{proof}

\begin{rem} \label{rem1.4}
Under the condition of Theorem \ref{thm1.3}, according to \cite{BBD}, Proposition 1.3.3,
there exists a functor $(-)^{\geq n}: \mcal{T} \arr \mcal{T}^{\geq n}(C)$
(resp., $(-)^{\leq n}: \mcal{T} \arr \mcal{T}^{\leq n}(C)$) which is the right
(resp., left) adjoint of the natural embedding functor
$\mcal{T}^{\geq n}(C) \arr \mcal{T}$ 
(resp., $\mcal{T}^{\leq n}(C) \arr \mcal{T}$).
\end{rem}

For an object $C$ in a triangulated category $\mcal{T}$ and integers $s \leq t$, let
$\mcal{T}^{[s]}(C) = \mcal{T}^{0}(C)[-s]$,
$\mcal{T}^{[s,t]}(C) = \mcal{T}^{\leq t}(C) \cap \mcal{T}^{\geq s}(C)$,  and
$\mcal{T}^{\mrm{b}}(C) = ({\bigcup}_{n \in \mbb{Z}}\mcal{T}^{\leq n}(C))
\cap \\ ({\bigcup}_{n \in \mbb{Z}}\mcal{T}^{\geq n}(C))$.
An object $M$ of an abelian category $\mcal{A}$ is called small provided that
$\opn{Hom}_{\mcal{A}}(M,  -)$ commutes with direct sums in $\mcal{A}$.

\begin{cor} \label{cor1.5}
Let $\mcal{A}$ be an abelian category satisfying the condition Ab4 
(i.e. direct sums of exact sequences are exact)
and $T^{\centerdot}$ a bounded  complex of small projective objects of $\mcal{A}$ satisfying
\begin{enumerate}
\rmitem{i} $\{T^{\centerdot}[i] : i \in \mbb{Z}\}$ is a generating set
for $\cat{D}(\mcal{A})$,
\rmitem{ii} $\opn{Hom}_{\cat{D}(\mcal{A})}(T^{\centerdot}, T^{\centerdot}[i]) = 0$
for $i \not= 0$.
\end{enumerate}
If either of the following conditions (1) or (2) is satisfied, then we have an equivalence
of triangulated categories
\[
\cat{D}(\mcal{A})^{\mrm{b}}(T^{\centerdot}) \cong
\cat{D}^{\mrm{b}}(\cat{Mod}B),
\]
where $B = \opn{End}_{\cat{D}(\mcal{A})}(T^{\centerdot})^{\mrm{op}}$.
\begin{enumerate}
\rmitem{1} $\mcal{A}$ has enough projectives.
\rmitem{2} $\mcal{A}$ has enough injectives and 
$\cat{D}(\mcal{A})^{\geq 0}(T^{\centerdot}) \subset \cat{D}^{+}(\mcal{A})$.
\end{enumerate}
Moreover, if $\cat{D}(\mcal{A})^{0}(T^{\centerdot}) \subset \cat{D}^{\mrm{b}}(\mcal{A})$,
then we have an equivalence 
\[
\cat{D}^{\mrm{b}}(\mcal{A}) \cong 
\cat{D}^{\mrm{b}}(\cat{Mod}B).
\]
\end{cor}

\begin{proof}
According to \cite{BN}, Corollary 1.7, $\cat{D}(\mcal{A})$ contains direct sums.
Since $T^{\centerdot}$ is a bounded  complex of small projective objects of $\mcal{A}$,
$T^{\centerdot}$ is a compact object in $\cat{D}(\mcal{A})$.  By Theorem \ref{thm1.3}
$\cat{D}(\mcal{A})$ has a $t$-structure $(\cat{D}(\mcal{A})^{\leq 0}(T^{\centerdot}), 
\cat{D}(\mcal{A})^{\geq 0}(T^{\centerdot}))$, and 
$\opn{Hom}_{\cat{D}(\mcal{A})}(T^{\centerdot}, -) : \cat{D}(\mcal{A})^{0}(T^{\centerdot})
\arr \cat{Mod}B$ is an equivalence.

(1) By the construction of $X^{\geq r}$ in Proposition \ref{prop1.1},
$\cat{D}^{-}(\mcal{A})$ also has a $t$-structure
$(\cat{D}^{-}(\mcal{A})^{\leq 0}(T^{\centerdot}),  
\cat{D}^{-}(\mcal{A})^{\geq 0}(T^{\centerdot}))$
and hence by Theorem \ref{thm1.3} (3) we have
$\cat{D}^{-}(\mcal{A})^{0}(T^{\centerdot}) = 
\cat{D}(\mcal{A})^{0}(T^{\centerdot})$.
According to \cite{Ri}, Proposition 10.1, we have a fully faithful $\partial$-functor
$F' : \cat{D}^{-}(\cat{Mod}B) \arr \cat{D}^{-}(\mcal{A})$.
Also, since $F'(B) = T^{\centerdot}$, $F'$ sends $B$-modules to objects in 
$\cat{D}(\mcal{A})^{0}(T^{\centerdot})$.   Then we have a fully faithful $\partial$-functor
\[
F : \cat{D}^{\mrm{b}}(\cat{Mod}B) \arr \cat{D}(\mcal{A}) ,
\]
which sends $B$-modules to objects in $\cat{D}(\mcal{A})^{0}(T^{\centerdot})$. 
For any $X \in \cat{D}(\mcal{A})^{\mrm{b}}(T^{\centerdot})$, there exist integers $m \leq n$
such that $X \in \cat{D}(\mcal{A})^{[m,n]}(T^{\centerdot})$.
Let $ l = n - m$.  If $l = 0$, then there exist obviously a $B$-module $M$
and an integer $s$ such that $X \cong F(M[s])$.
If $l > 0$, then we have a distinguished triangle
\[
X^{\leq n-1} \arr X \arr X^{\geq n} \arr
\]
with $X^{\geq n} \in \cat{D}(\mcal{A})^{[n]}(T^{\centerdot})$ and
$X^{\leq n-1} \in \cat{D}(\mcal{A})^{[m,n-1]}(T^{\centerdot})$.
Since $F$ is full, by induction on $l$, there exists
$U^{\centerdot} \in \cat{D}^{\mrm{b}}(\cat{Mod}B)$ such that $X \cong 
F(U^{\centerdot})$.

(2) By the assumption,
$\cat{D}^{+}(\mcal{A})$ also has a $t$-structure
$(\cat{D}^{+}(\mcal{A})^{\leq 0}(T^{\centerdot}),  \\
\cat{D}^{+}(\mcal{A})^{\geq 0}(T^{\centerdot}))$.  Thus 
$\cat{D}^{+}(\mcal{A})^{\mrm{b}}(T^{\centerdot}) =
\cat{D}(\mcal{A})^{\mrm{b}}(T^{\centerdot})$, and hence
$\cat{D}^{+}(\mcal{A})^{0}(T^{\centerdot}) = \\
\cat{D}(\mcal{A})^{0}(T^{\centerdot})$.
By \cite{BBD}, Section 3, we have a $\partial$-functor
$\opn{real} : \cat{D}^{\mrm{b}}(\cat{D}(\mcal{A})^{0}(T^{\centerdot})) 
\arr \cat{D}^{+}(\mcal{A})$, and then we have a $\partial$-functor
\[
F : \cat{D}^{\mrm{b}}(\cat{Mod}B) \arr \cat{D}(\mcal{A}) ,
\]
which sends $B$-modules to objects in $\cat{D}(\mcal{A})^{0}(T^{\centerdot})$. 
Let $f \in \opn{Hom}_{\cat{D}(\mcal{A})}(X^{\centerdot}, Y^{\centerdot}[n])$
with $X^{\centerdot}, Y^{\centerdot} \in \cat{D}(\mcal{A})^{0}(T^{\centerdot})$ and 
$n > 0$.
Take a distinguished triangle in $\cat{D}^{+}(\mcal{A})$
\[
X^{\centerdot}_1 \arr V^{\centerdot} \xarr{t} X^{\centerdot} \arr
\]
such that $V^{\centerdot}$ is a direct sum of copies of $T^{\centerdot}$ and
$\opn{Hom}_{\cat{D}(\mcal{A})}(T^{\centerdot},t)$ is an epimorphism.
By easy calculation, $X^{\centerdot}_1 \in \cat{D}(\mcal{A})^{0}(T^{\centerdot})$, and hence
we get an exact
sequence in $\cat{D}(\mcal{A})^{0}(T^{\centerdot})$
\[
0 \arr X^{\centerdot}_1 \arr V^{\centerdot} \xarr{t} X^{\centerdot} \arr 0 .
\]
Since $\opn{Hom}_{\cat{D}(\mcal{A})}(T^{\centerdot},Y^{\centerdot}[n]) = 0$, we have $ft = 0$, 
i.e. $t$ effaces $f$.
Thus the epimorphic version of effacibility in \cite{BBD}, Proposition 3.1.16
can be applied.

Finally, it is easy to see that $\cat{D}(\mcal{A})^{0}(T^{\centerdot}) \subset \cat{D}^{\mrm{b}}(\mcal{A})$
implies $\cat{D}^{\mrm{b}}(\mcal{A}) = 
\cat{D}(\mcal{A})^{\mrm{b}}(T^{\centerdot})$.
\end{proof}

\section{Torsion Theories for Abelian Categories} \label{s2}

Throughout this section,  we fix the following notation. 
Let $\mcal{A}$ be an abelian category satisfying the condition Ab4, and let
$d_{P}^{- 1} : P^{- 1} \arr P^{0}$ 
be a morphism in $\mcal{A}$ with the $P^{i}$ being small projective objects 
of $\mcal{A}$, and denote by $P^{\centerdot}$ the mapping cone of $d_{P}^{- 1}$.
We set $\mcal{C}(P^{\centerdot}) = \cat{D}(\mcal{A})^{0}(P^{\centerdot})$,
$B = \opn{Hom}_{\cat{D}(\mcal{A})}(P^{\centerdot})^{\mrm{op}}$, and
define a pair of full subcategories of $\mcal{A}$
\[\begin{aligned}
\mcal{X}(P^{\centerdot}) & = \{ X \in \mcal{A} :
\opn{Hom}_{\cat{D}(\mcal{A})}(P^{\centerdot}, X[1]) = 0 \} , \\
\mcal{Y}(P^{\centerdot}) & = \{ X \in \mcal{A} :
\opn{Hom}_{\cat{D}(\mcal{A})}(P^{\centerdot}, X) = 0 \}.
\end{aligned} \]
For any $X \in \mcal{A}$, we define a subobject of $X$
\[
{\tau}(X) = {\sum}_{f \in 
\opn{Hom}_{\mcal{A}}(\opn{H}^{0}(P^{\centerdot}),X)}\opn{Im}f
\]
and an exact sequence in $\mcal{A}$
\[
(e_X) : 0 \arr {\tau}(X) \xarr{j_{X}} X \arr {\pi}(X) \arr 0.
\]

\begin{rem} \label{rem2.1}
It is easy to see that $P^{\centerdot}$ is a compact object of $\cat{D}(\mcal{A})$, and
we have $\mcal{X}(P^{\centerdot}) = \cat{D}(\mcal{A})^{\leq 0}(P^{\centerdot})
\cap\mcal{A}$ and
$\mcal{Y}(P^{\centerdot}) = \cat{D}(\mcal{A})^{\geq 1}(P^{\centerdot})
\cap\mcal{A}$.
\end{rem}

\begin{lem} \label{lem2.2}
For any $X \in \mcal{A}$, the following hold.
\begin{enumerate}
\rmitem{1} $\opn{Ker}(\opn{Hom}_{\mcal{A}}(d_{P}^{- 1},X)) \cong
\opn{Hom}_{\cat{D}(\mcal{A})}(P^{\centerdot},X)$.
\rmitem{2} $\opn{Cok}(\opn{Hom}_{\mcal{A}}(d_{P}^{- 1},X)) \cong
\opn{Hom}_{\cat{D}(\mcal{A})}(P^{\centerdot},X[1])$.
\end{enumerate}
\end{lem}

\begin{lem} \label{lem2.3}
For any $X \in \mcal{A}$, the following hold.
\begin{enumerate}
\rmitem{1} $\opn{Hom}_{\cat{D}(\mcal{A})}(P^{\centerdot},X[n]) = 0$
for $n >1$ and $n < 0$.
\rmitem{2} $\opn{Hom}_{\cat{D}(\mcal{A})}(P^{\centerdot},X) \cong
\opn{Hom}_{\mcal{A}}(\opn{H}^{0}(P^{\centerdot}),X)$.
\end{enumerate}
\end{lem}

\begin{lem} \label{lem2.4}
The following hold.
\begin{enumerate}
\rmitem{1}  $\mcal{X}(P^{\centerdot})$ is closed under factor objects and direct sums.
\rmitem{2}  $\mcal{Y}(P^{\centerdot})$ is closed under subobjects.
\rmitem{3}  For any $X \in \mcal{A}$, $\opn{Hom}_{A}(\opn{H}^{0}(P^{\centerdot}), j_X)$
is an isomorphism.
\end{enumerate}
\end{lem}

\begin{lem} \label{lem2.5}
For any $X^{\centerdot} \in \cat{D}(\mcal{A})$ and $n \in \mbb{Z}$, we have a 
functorial exact sequence

{\small \[
0 \arr \opn{Hom}_{\cat{D}(\mcal{A})}(P^{\centerdot}, 
\opn{H}^{n-1}(X^{\centerdot})[1]) \arr 
\opn{Hom}_{\cat{D}(\mcal{A})}(P^{\centerdot}, X^{\centerdot}[n]) \arr 
\opn{Hom}_{\cat{D}(\mcal{A})}(P^{\centerdot}, \opn{H}^{n}(X^{\centerdot})) \arr 0.
\]} \par\noindent
Moreover, the above short exact sequence commutes with direct sums.
\end{lem}

\begin{proof}
For $X^{\centerdot}[n] \in \cat{D}(\mcal{A})$, applying 
$\opn{Hom}_{\cat{D}(\mcal{A})}(-, X^{\centerdot}[n])$ to a distinguished triangle
\[
P^{- 1} \xarr{d_{P}^{- 1}} P^{0} \arr P^{\centerdot} \arr ,
\]
we have a short exact sequence
\[\begin{aligned}
0 \arr \opn{Cok}(\opn{Hom}_{\cat{D}(\mcal{A})}(d_{P}^{-1}, 
X^{\centerdot}[n- 1])) & \arr
\opn{Hom}_{\cat{D}(\mcal{A})}(P^{\centerdot}, X^{\centerdot}[n]) \\ 
& \arr
\opn{Ker}(\opn{Hom}_{\cat{D}(\mcal{A})}(d_{P}^{-1}, X^{\centerdot}[n]))
\arr 0.
\end{aligned}\]
Also, by Lemma \ref{lem2.2} we get
\[ \begin{aligned}
\opn{Ker}(\opn{Hom}_{\cat{D}(\mcal{A})}(d_{P}^{-1}, X^{\centerdot}[n])) 
& \cong \opn{Ker}(\opn{Hom}_{\mcal{A}}(d_{P}^{-1}, \opn{H}^{n}(X^{\centerdot}))) \\
& \cong \opn{Hom}_{\cat{D}(\mcal{A})}(P^{\centerdot}, \opn{H}^{n}(X^{\centerdot})),
\end{aligned}
\] \[
\begin{aligned}
\opn{Cok}(\opn{Hom}_{\cat{D}(\mcal{A})}(d_{P}^{-1},X^{\centerdot}[n-1])) 
& \cong \opn{Cok}(\opn{Hom}_{\mcal{A}}(d_{P}^{-1},\opn{H}^{n-1}(X^{\centerdot}))) \\
& \cong \opn{Hom}_{\cat{D}(\mcal{A})}(P^{\centerdot}, 
\opn{H}^{n-1}(X^{\centerdot})[1]).
\end{aligned} \]
Since the $P^{i}$ are small objects, the above short exact sequence 
commutes with direct sums.
\end{proof}

\begin{lem} \label{lem2.6}
The following are equivalent.
\begin{enumerate}
\rmitem{1} $\{P^{\centerdot}[i] : i \in \mbb{Z}\}$ is a generating set
for $\cat{D}(\mcal{A})$.
\rmitem{2} $\mcal{X}(P^{\centerdot})\cap\mcal{Y}(P^{\centerdot}) = \{ 0\}$.
\end{enumerate}
\end{lem}

\begin{proof} (1) $\Rightarrow$ (2).  For any $X \in \mcal{X}(P^{\centerdot}) 
\cap \mcal{Y}(P^{\centerdot})$,
by Lemma \ref{lem2.3} (1), $\opn{Hom}_{\cat{D}(\mcal{A})}(P^{\centerdot}, \\ X[n]) = 0$ 
for all $n \in \mbb{Z}$ and hence $X = 0$.

(2) $\Rightarrow$ (1).  Let $X^{\centerdot} \in \cat{D}(\mcal{A})$ with 
$\opn{Hom}_{\cat{D}(\mcal{A})}(P^{\centerdot}, X^{\centerdot}[n]) = 0$ for all 
$n \in \mbb{Z}$.  Then by Lemma \ref{lem2.5}, $\opn{H}^{n}(X^{\centerdot}) 
\in \mcal{X}(P^{\centerdot})\cap\mcal{Y}(P^{\centerdot}) = \{ 0\}$.
\end{proof}

\begin{lem} \label{lem2.7}
The following hold.
\begin{enumerate}
\rmitem{1} $\opn{H}^{0}(P^{\centerdot}) \in \mcal{X}(P^{\centerdot})$
if and only if $\opn{Hom}_{\cat{D}(\mcal{A})}(P^{\centerdot},P^{\centerdot}[i])
= 0$ for all $i > 0$.
\rmitem{2} $\opn{H}^{-1}(P^{\centerdot}) \in \mcal{Y}(P^{\centerdot})$
if and only if $\opn{Hom}_{\cat{D}(\mcal{A})}(P^{\centerdot},P^{\centerdot}[i])
= 0$ for all $i < 0$.
\end{enumerate}
\end{lem}

\begin{proof}
By Lemma \ref{lem2.5}.
\end{proof}

\begin{defn} \label{dfn2.8}
A pair $(\mcal{X}, \mcal{Y})$ of full subcategories $\mcal{X},
\mcal{Y}$ in an abelian category $\mcal{A}$ is called a torsion 
theory for $\mcal{A}$ provided that the following conditions are satisfied
(see e.g. \cite{Di} for details):
\begin{enumerate}
\rmitem{i} $\mcal{X} \cap \mcal{Y} = \{ 0\}$;
\rmitem{ii} $\mcal{X}$ is closed under factor objects;
\rmitem{iii} $\mcal{Y}$ is closed under subobjects;
\rmitem{iv} for any object $X$ of $\mcal{A}$, there exists an exact 
sequence $0 \arr X' \arr X \arr X'' \arr 0$ in $\mcal{A}$ with 
$X' \in \mcal{X}$ and $X'' \in \mcal{Y}$.
\end{enumerate}
\end{defn}

\begin{rem} \label{rem2.9}
Let $\mcal{A}$ be an abelian category and $(\mcal{X}, \mcal{Y})$ a torsion 
theory for $\mcal{A}$.  Then for any $Z \in \mcal{A}$, the following hold.
\begin{enumerate}
\rmitem{1} $Z \in \mcal{X}$ if and only if 
$\opn{Hom}_{\mcal{A}}(Z, \mcal{Y}) = 0$.
\rmitem{2}  $Z \in \mcal{Y}$ if and only if
$\opn{Hom}_{\mcal{A}}(\mcal{X}, Z) = 0$.
\end{enumerate}
\end{rem}

\begin{thm} \label{thm2.10}
The following are equivalent.
\begin{enumerate}
\rmitem{1} $\{P^{\centerdot}[i] : i \in \mbb{Z}\}$ is a generating set
for $\cat{D}(\mcal{A})$ and $\opn{Hom}_{\cat{D}(\mcal{A})}(P^{\centerdot},
P^{\centerdot}[i]) = 0$ for all $i > 0$.
\rmitem{2} $\mcal{X}(P^{\centerdot}) \cap \mcal{Y}(P^{\centerdot}) = \{ 0\}$
and $\opn{H}^{0}(P^{\centerdot}) \in \mcal{X}(P^{\centerdot})$.
\rmitem{3} $\mcal{X}(P^{\centerdot}) \cap \mcal{Y}(P^{\centerdot}) = \{ 0\}$
and $\tau(X) \in \mcal{X}(P^{\centerdot})$, $\pi(X) \in 
\mcal{Y}(P^{\centerdot})$ for all $X \in \mcal{A}$.
\rmitem{4} $(\mcal{X}(P^{\centerdot}), \mcal{Y}(P^{\centerdot}))$ is a torsion 
theory for $\mcal{A}$.
\end{enumerate}
\end{thm}

\begin{proof} 
(1) $\Leftrightarrow$ (2). By Lemmas \ref{lem2.6} and \ref{lem2.7} (1).

(2) $\Rightarrow$ (3).
Let $X \in \mcal{A}$.  Since $\opn{H}^{0}(P^{\centerdot})
\in \mcal{X}(P^{\centerdot})$, it follows that
$\tau(X) \in \mcal{X}(P^{\centerdot})$.
Next, apply $\opn{Hom}_{\cat{D}(\mcal{A})}(P^{\centerdot}, -)$ to the
canonical exact sequence $(e_X)$.
It then follows by Lemmas \ref{lem2.3} (2) and \ref{lem2.4} (3)
that $\opn{Hom}_{\cat{D}(\mcal{A})}(P^{\centerdot}, j_{X})$ is an isomorphism.
Thus $\opn{Hom}_{\cat{D}(\mcal{A})}(P^{\centerdot}, \pi(X)) = 0$ and hence
$\pi(X) \in \mcal{Y}(P^{\centerdot})$.

(3) $\Rightarrow$ (4).  Obvious.

(4) $\Rightarrow$ (1).  By Lemmas \ref{lem2.3} (2), \ref{lem2.6}, \ref{lem2.7} (1)
and Remark \ref{rem2.9} (1).
\end{proof}

\begin{defn} \label{dfn2.11}
For a complex $X^{\centerdot} = (X^{i}, d^{i})$, we define 
the following truncations:
\[
\begin{aligned}
{\sigma}_{> n}(X^{\centerdot}) & : \ldots \arr 0 \arr \opn{Im}d^{n} 
\arr X^{n+1} \arr X^{n+2} \arr \ldots ,\\
{\sigma}_{\leq n}(X^{\centerdot}) & : \ldots \arr X^{n-2} \arr X^{n-1}
\arr \opn{Ker}d^{n} \arr 0 \arr \ldots ,\\
{\sigma '}_{\geq n}(X^{\centerdot}) & : \ldots \arr 0 \arr \opn{Cok}d^{n-1}
\arr X^{n+1}\arr X^{n+2}\arr \ldots ,\\
{\sigma '}_{< n}(X^{\centerdot}) & : \ldots \arr X^{n-2} \arr X^{n-1}
\arr \opn{Im}d^{n-1} \arr 0 \arr \ldots .
\end{aligned}
\]
\end{defn}

\begin{lem} \label{lem2.12}
For any $X^{\centerdot} \in \cat{D}(\mcal{A})$ with 
$\opn{H}^{n}(X^{\centerdot}) = 0$ for $n > 0$ and $n < - 1$, there exists a 
distinguished triangle in $\cat{D}(\mcal{A})$ of the form
\[
\opn{H}^{-1}(X^{\centerdot})[1] \arr X^{\centerdot} \arr 
\opn{H}^{0}(X^{\centerdot}) \arr .
\]
\end{lem}

\begin{proof}
We have exact sequences in $\cat{C}(\mcal{A})$
\[     
0 \arr \sigma_{\leq -1}(X^{\centerdot}) \arr X^{\centerdot} \arr 
\sigma_{> -1}(X^{\centerdot}) \arr 0,
\] \[
0 \arr {\sigma '}_{< 0}(\sigma_{> -1}(X^{\centerdot})) \arr 
\sigma_{> -1}(X^{\centerdot}) \arr {\sigma '}_{\geq 0}(X^{\centerdot}) \arr 0.
\]
Also, $\sigma_{\leq -1}(X^{\centerdot}) \cong \opn{H}^{-1}(X^{\centerdot})[1]$,
${\sigma '}_{< 0}(\sigma_{> -1}(X^{\centerdot})) \cong 0$ and 
${\sigma '}_{\geq 0}(X^{\centerdot}) \cong \opn{H}^{0}(X^{\centerdot})$
in $\cat{D}(\mcal{A})$.  Thus we get a desired 
distinguished triangle in $\cat{D}(\mcal{A})$.
\end{proof}

\begin{lem} \label{lem2.13}
Assume $\mcal{X}(P^{\centerdot}) \cap \mcal{Y}(P^{\centerdot}) 
= \{ 0 \}$.  Then for any $X^{\centerdot} \in \cat{D}(\mcal{A})$,
the following are equivalent.
\begin{enumerate}
\rmitem{1} $X^{\centerdot} \in \mcal{C}(P^{\centerdot})$.
\rmitem{2} $\opn{H}^n(X^{\centerdot}) = 0$ for $n > 0$ and $n < - 1$, 
$\opn{H}^{0}(X^{\centerdot}) \in \mcal{X}(P^{\centerdot})$ and 
$\opn{H}^{-1}(X^{\centerdot}) \in \mcal{Y}(P^{\centerdot})$.
\end{enumerate}
\end{lem}

\begin{proof}
By Lemma \ref{lem2.5}.
\end{proof}

\begin{rem} \label{rem2.14}
Let $\mcal{A}$ be an abelian category and $\mcal{X}, \mcal{Y}$ full 
subcategories of $\mcal{A}$.  Then the pair $(\mcal{X}, \mcal{Y})$ 
is a torsion theory for $\mcal{A}$ 
if and only if the following two conditions are satisfied:
\begin{enumerate}
\rmitem{i} $\opn{Hom}_{\mcal{A}}(\mcal{X}, \mcal{Y}) = 0$;
\rmitem{ii} for any object $X$ in $\mcal{A}$, there exists an exact sequence 
$0 \arr X' \arr X \arr X'' \arr 0$ in $\mcal{A}$ with $X' \in \mcal{X}$ 
and $X'' \in \mcal{Y}$.
\end{enumerate}
\end{rem}

\begin{thm} \label{thm2.15}
Assume $\mcal{X}(P^{\centerdot}) \cap
\mcal{Y}(P^{\centerdot}) = \{ 0 \}$ and $\opn{H}^{0}(P^{\centerdot}) \in 
\mcal{X}(P^{\centerdot})$.  Then the following hold.
\begin{enumerate}
\rmitem{1} $\mcal{C}(P^{\centerdot})$ is admissible abelian.
\rmitem{2} The functor
\[
\opn{Hom}_{\cat{D}(\mcal{A})}(P^{\centerdot}, -) : 
\mcal{C}(P^{\centerdot}) \arr \cat{Mod}B
\]
is an  equivalence.
\rmitem{3} $(\mcal{Y}(P^{\centerdot})[1], \mcal{X}(P^{\centerdot}))$ is a torsion 
theory for $\mcal{C}(P^{\centerdot})$.
\end{enumerate}
\end{thm}

\begin{proof}
(1) and (2) According to Theorem \ref{thm2.10},  Theorem \ref{thm1.3}
can be applied.

(3) Note first that by Lemma \ref{lem2.13} we have 
$\mcal{X}(P^{\centerdot}) \subset \mcal{C}(P^{\centerdot})$ and 
$\mcal{Y}(P^{\centerdot})[1] \subset \mcal{C}(P^{\centerdot})$.  Also, 
it is trivial that $\opn{Hom}_{\cat{D}(\mcal{A})}(\mcal{Y}(P^{\centerdot})[1], 
\mcal{X}(P^{\centerdot})) = 0$.  Let 
$X^{\centerdot} \in \mcal{C}(P^{\centerdot})$.  Then by Lemmas \ref{lem2.12} 
and \ref{lem2.13} we have a distinguished triangle in 
$\cat{D}(\mcal{A})$ of the form
\[
\opn{H}^{-1}(X^{\centerdot})[1] \arr X^{\centerdot} \arr 
\opn{H}^{0}(X^{\centerdot}) \arr .
\]
It follows that the sequence in $\mcal{C}(P^{\centerdot})$
\[
0 \arr \opn{H}^{-1}(X^{\centerdot})[1] \arr X^{\centerdot} \arr 
\opn{H}^{0}(X^{\centerdot}) \arr 0
\]
is exact.  Thus by Remark \ref{rem2.14} $(\mcal{Y}(P^{\centerdot})[1], 
\mcal{X}(P^{\centerdot}))$ is a torsion theory for $\mcal{C}(P^{\centerdot})$.
\end{proof}

\begin{prop} \label{prop2.16}
Assume $P^{\centerdot}$ satisfies the conditions
\begin{enumerate}
\rmitem{i} $\{P^{\centerdot}[i]: i \in \mbb{Z}\}$ is a generating set for 
$\cat{D}(\mcal{A})$,
\rmitem{ii} $\opn{Hom}_{\cat{D}(\mcal{A})}(P^{\centerdot}, P^{\centerdot}[i]) = 0$
for $i \not= 0$.
\end{enumerate}
If $\mcal{A}$ has either enough projectives or enough injectives,
then we have an equivalence of triangulated categories 
\[
\cat{D}^{\mrm{b}}(\mcal{A}) \cong \cat{D}^{\mrm{b}}(\cat{Mod}B)  .
\]
\end{prop} 

\begin{proof}
Let $X^{\centerdot} \in \cat{D}(\mcal{A})$.
According to Lemma \ref{lem2.5} and Theorem \ref{thm2.10}, it is easy
to see that if $X^{\centerdot}$ belongs to $\cat{D}(\mcal{A})^{\geq 0}(P^{\centerdot})$ 
(resp., $\mcal{C}(P^{\centerdot})$), then $\opn{H}^{n}(X^{\centerdot}) = 0$
for $n < -1$ (resp., $n < -1$ and $n > 0$).
Thus we have
\[
\cat{D}(\mcal{A})^{\geq 0}(P^{\centerdot}) \subset \cat{D}^{+}(\mcal{A})
\enskip\text{and}\enskip
\mcal{C}(P^{\centerdot}) \subset \cat{D}^{\mrm{b}}(\mcal{A}) ,
\]
so that Corollary \ref{cor1.5} can be applied.
\end{proof}

\section{Torsion Theories for Module Categories} \label{s3}

In this section, we apply results of Section \ref{s2} to the case  of module categories.
In and after this section,  $R$  is a commutative ring and $I$ is an 
injective cogenerator in the category of $R$-modules.  We set 
$D = \opn{Hom}_R(-, I)$.  Let 
$A$ be an $R$-algebra and denote by $\cat{Proj}A$ (resp., $\cat{proj}A$) 
the full additive subcategory of 
$\cat{Mod}A$ consisting of projective (resp., finitely generated projective) 
modules.  We denote by $A^{\mrm{op}}$ the opposite ring of 
$A$ and consider right $A$-modules as left $A^{\mrm{op}}$-modules.
Also, we denote by $(-)^{*}$ both the $A$-dual functors $\opn{Hom}_{A}(-, A)$ 
and set 
$\nu = D\circ (-)^{*}$.
 
It is well known that, in a module category, the small projective objects are just 
the finitely generated projective modules.  In the following, we deal with the case where 
$\mcal{A} = \cat{Mod}A$ and use the same notation as in Section \ref{s2}.

\begin{lem} \label{lem3.1}
For any $X \in \cat{Mod}A$, we have
\[
\opn{Hom}_{\cat{D}(\cat{Mod}A)}(P^{\centerdot}, X[1]) \cong 
\opn{H}^{1}((P^{\centerdot})^{*}){\otimes}_{A}X .
\]
\end{lem}

\begin{proof}  We have 
\[\begin{aligned}
\opn{Hom}_{\cat{D}(\cat{Mod}A)}(P^{\centerdot}, X[1]) & \cong
\opn{Hom}_{\cat{K}(\cat{Mod}A)}(P^{\centerdot}, X[1]) \\
& \cong \opn{H}^{1}(\opn{Hom}^{\centerdot}_{A}(P^{\centerdot}, X)) \\
& \cong \opn{H}^{1}((P^{\centerdot})^{*} {\otimes}^{\centerdot}_{A} X) \\
& \cong \opn{H}^{1}((P^{\centerdot})^{*}) {\otimes}_{A} X.
\end{aligned} \]
\end{proof}

\begin{lem} \label{lem3.2}
The following hold.
\begin{enumerate}
\rmitem{1} $\mcal{X}(P^{\centerdot}) = 
\opn{Ker}(\opn{H}^1((P^{\centerdot})^{*}) {\otimes}_A -)$.
\rmitem{2} $\mcal{Y}(P^{\centerdot}) = 
\opn{Ker}(\opn{Hom}_{A}(\opn{H}^0(P^{\centerdot}), -))$.
\end{enumerate}
\end{lem}

\begin{proof}
By Lemmas \ref{lem2.3} (2) and \ref{lem3.1}.
\end{proof}

\begin{lem} \label{lem3.3}
The following hold.
\begin{enumerate}
\rmitem{1} $D(\opn{H}^{1}((P^{\centerdot})^{*})) \cong \opn{H}^{-1}(\nu(P^{\centerdot}))$.
\rmitem{2} $\mcal{X}(P^{\centerdot}) = \opn{Ker}(\opn{Hom}_{A}(-, \opn{H}^{-1}(\nu(P^{\centerdot}))))$
and hence 
$\opn{H}^{0}(P^{\centerdot}) \in \mcal{X}(P^{\centerdot})$
if and only if
$\opn{H}^{-1}(\nu(P^{\centerdot})) \in \mcal{Y}(P^{\centerdot})$.
\rmitem{3} $\opn{Ker}(\opn{Tor}^{A}_{1}(\opn{H}^{1}((P^{\centerdot})^{*}), -)) = 
\opn{Ker}(\opn{Ext}^{1}_{A}(-, \opn{H}^{-1}(\nu(P^{\centerdot}))))$.
\end{enumerate}
\end{lem}
     
\begin{proof}
We have $D(\opn{H}^{1}((P^{\centerdot})^{*})) \cong 
\opn{H}^{-1}(D((P^{\centerdot})^{*})) = 
\opn{H}^{-1}(\nu(P^{\centerdot}))$ and for any $X \in \cat{Mod}A$ we have 
\[
\begin{aligned}
D(\opn{H}^{1}((P^{\centerdot})^{*}) {\otimes}_{A} X) & \cong \opn{Hom}_{A}(X, 
\opn{H}^{-1}(\nu(P^{\centerdot}))), \\
D(\opn{Tor}^{A}_{1}((\opn{H}^{1}((P^{\centerdot})^{*}), X))) & 
\cong \opn{Ext}^{1}_{A}(X, \opn{H}^{-1}(\nu(P^{\centerdot}))).
\end{aligned} \]
\end{proof}

\begin{lem} \label{lem3.4}
The following hold.
\begin{enumerate}
\rmitem{1} $\mcal{X}(P^{\centerdot}) \subset 
\opn{Ker}(\opn{Ext}^{1}_{A}(\opn{H}^{0}(P^{\centerdot}), -))$.
\rmitem{2} $\mcal{Y}(P^{\centerdot}) \subset 
\opn{Ker}(\opn{Tor}^{A}_{1}(\opn{H}^{1}((P^{\centerdot})^{*}), -))$.
\end{enumerate}
\end{lem}

\begin{proof} This is due essentially to Auslander \cite{Au}.
We have an exact sequence in $\cat{Mod}A$
\[
0 \arr \opn{H}^{-1}(P^{\centerdot}) \arr P^{- 1} \arr P^{0} \arr 
\opn{H}^{0}(P^{\centerdot}) \arr 0
\]
with the $P^{i}$ finitely generated projective, and an exact 
sequence in $\cat{Mod}A^{\mrm{op}}$
\[
0 \arr \opn{H}^{0}(P^{\centerdot})^{*} \arr P^{0 *} \arr P^{-1 *} \arr 
\opn{H}^{1}((P^{\centerdot})^{*}) \arr 0
\]
with the $P^{i *}$ finitely generated projective.

(1) Let $X \in \cat{Mod}A$.  For any $M \in \cat{Mod}A^{\mrm{op}}$, 
we have a functorial homomorphism
\[
{\theta}_{M} : M{\otimes}_{A}X \arr \opn{Hom}_{A}(M^{*}, X),
m{\otimes}x \mapsto (h \mapsto h(m)x)
\]
which is an isomorphism if $M$ is finitely generated projective.
Since the $P^{i}$ are reflexive, we have $\opn{H}^{0}(P^{\centerdot}) \cong 
\opn{H}^{0}((P^{\centerdot})^{* *})$ and $\opn{H}^{-1}(P^{\centerdot}) \cong 
\opn{H}^{1}((P^{\centerdot})^{*})^{*}$.  We have a commutative diagram
\[ \begin{CD}
P^{0 *} {\otimes}_{A} X @>>> P^{-1 *} {\otimes}_{A} X @>>>
\opn{H}^{1}((P^{\centerdot})^{*}) {\otimes}_{A} X @>>> 0 \\
@V{\theta}_{P^{0 *}}VV @VV{\theta}_{P^{-1 *}}V @VVV \\
\opn{Hom}_{A}(P^{0 * *}, X) @>>> \opn{Hom}_{A}(P^{-1 * *}, X) @>>>
\opn{Hom}_{A}(\opn{H}^{-1}(P^{\centerdot}), X) 
\end{CD} \]
with the top row exact.  Since the ${\theta}_{P^{i *}}$ are isomorphisms, 
$\opn{Ext}_{A}^{1}(\opn{H}^{0}(P^{\centerdot}), X)$ is embedded in 
$\opn{H}^{1}((P^{\centerdot})^{*}){\otimes}_{A}X$.  
The assertion follows by Lemma \ref{lem3.2}.

(2) Let $X \in \cat{Mod}A$.  For any $Y \in \cat{Mod}A$, we have a 
functorial homomorphism
\[
{\eta}_{Y} : Y^{*} {\otimes}_{A} X \arr \opn{Hom}_{A}(Y, X),
h {\otimes} x \mapsto (y \mapsto h(y)x)
\]
which is an isomorphism if $Y$ is finitely generated projective.
We have a commutative diagram
\[ \begin{CD}
@. \opn{H}^{0}(P^{\centerdot})^{*} {\otimes}_{A} X @>>> P^{0 *} {\otimes}_{A} X
@>>> P^{- 1 *} {\otimes}_{A} X \\
@. @VVV @V{\eta}_{P^{0}}VV @VV{\eta}_{P^{-1}}V \\ 
0 @>>> \opn{Hom}_{A}(\opn{H}^{0}(P^{\centerdot}), X) @>>> \opn{Hom}_{A}(P^{0}, X)
@>>> \opn{Hom}_{A}(P^{- 1}, X)
\end{CD} \]
with the bottom row exact.  Since the ${\eta}_{P^{i}}$ are isomorphisms, 
$\opn{Tor}^{A}_{1}(\opn{H}^{1}((P^{\centerdot})^{*}), X)$ is a homomorphic image of 
$\opn{Hom}_{A}(\opn{H}^{0}(P^{\centerdot}), X)$.  
The assertion follows by Lemma \ref{lem3.2}.
\end{proof}

\begin{thm} \label{thm3.5}
The following are equivalent.
\begin{enumerate}
\rmitem{1} $\mcal{X}(P^{\centerdot}) \cap \mcal{Y}(P^{\centerdot}) = \{ 0\}$
and $\opn{H}^{0}(P^{\centerdot}) \in \mcal{X}(P^{\centerdot})$.
\rmitem{2} $\mcal{X}(P^{\centerdot}) \cap \mcal{Y}(P^{\centerdot}) = \{ 0\}$
and $\tau(X) \in \mcal{X}(P^{\centerdot})$, $\pi(X) \in 
\mcal{Y}(P^{\centerdot})$ for all $X \in \cat{Mod}A$.
\rmitem{3} $(\mcal{X}(P^{\centerdot}), \mcal{Y}(P^{\centerdot}))$ is a torsion 
theory for $\cat{Mod}A$.
\rmitem{4} $\mcal{X}(P^{\centerdot})$ consists of the modules generated by 
$\opn{H}^{0}(P^{\centerdot})$ and $\mcal{Y}(P^{\centerdot})$ consists of the modules 
cogenerated by $\opn{H}^{-1}(\nu(P^{\centerdot}))$.
\end{enumerate}
\end{thm}

\begin{proof} (1) $\Leftrightarrow$ (2) $\Leftrightarrow$ (3).  By Theorem \ref{thm2.10}.

(3) $\Rightarrow$ (4).  Since $\opn{Hom}_{A}(\opn{H}^{0}(P^{\centerdot}), -)$
vanishes on $\mcal{Y}(P^{\centerdot})$, $\opn{H}^{0}(P^{\centerdot}) \in 
\mcal{X}(P^{\centerdot})$.
Thus $\mcal{X}(P^{\centerdot})$ contains the modules generated by 
$\opn{H}^{0}(P^{\centerdot})$.  Conversely, let $X \in \mcal{X}(P^{\centerdot})$.
Then, since (1) implies (2), $\pi(X) \in \mcal{Y}(P^{\centerdot})$ and hence
$\opn{Hom}_{A}(X, \pi(X)) = 0$. Thus $X = {\tau}(X)$, which is generated by 
$\opn{H}^{0}(P^{\centerdot})$.  Next, since by Lemma 
\ref{lem3.3} (2) $\opn{H}^{-1}(\nu(P^{\centerdot})) \\ \in \mcal{Y}(P^{\centerdot})$, 
$\mcal{Y}(P^{\centerdot})$ contains the modules cogenerated by 
$\opn{H}^{-1}(\nu(P^{\centerdot}))$.  Conversely, let $X \in 
\mcal{Y}(P^{\centerdot})$.  Take a set of generators 
$\{{f}_{\lambda}\}_{\lambda \in \Lambda}$ for an 
$R$-module $\opn{Hom}_{A}(X, \\ \opn{H}^{-1}(\nu(P^{\centerdot})))$ and set
\[     
f : X \arr \opn{H}^{-1}(\nu(P^{\centerdot}))^{\Lambda}, 
x \mapsto (f_{\lambda}(x))_{\lambda \in \Lambda}.
\]
It is obvious that $\opn{Hom}_{A}(f, \opn{H}^{-1}(\nu(P^{\centerdot})))$ is surjective.
Also, by Lemmas \ref{lem3.3} (3) and \ref{lem3.4}(2) we have
$\opn{Ext}^{1}_{A}(\opn{Im}f, \opn{H}^{-1}(\nu(P^{\centerdot}))) = 0$.  Applying 
$\opn{Hom}_{A}(-, \opn{H}^{-1}(\nu(P^{\centerdot})))$ to the canonical exact sequence
\[
0 \arr \opn{Ker}f \arr X \arr \opn{Im}f \arr 0,
\]
we get $\opn{Hom}_{A}(\opn{Ker}f, \opn{H}^{-1}(\nu(P^{\centerdot}))) = 0$.  Thus 
$\opn{Ker}f \in \mcal{X}(P^{\centerdot}) \cap \mcal{Y}(P^{\centerdot})$
and hence $\opn{Ker}f = 0$.

(4) $\Rightarrow$ (1).  By Lemma \ref{lem3.3} (2).
\end{proof}

\begin{cor} \label{cor3.6} 
The following are equivalent.
\begin{enumerate}
\rmitem{1} $P^{\centerdot}$ is a tilting complex.
\rmitem{2} $\mcal{X}(P^{\centerdot}) \cap \mcal{Y}(P^{\centerdot}) = \{ 0\}, 
\opn{H}^{0}(P^{\centerdot}) \in \mcal{X}(P^{\centerdot})$ and 
$\opn{H}^{-1}(P^{\centerdot}) \in \mcal{Y}(P^{\centerdot})$.
\rmitem{3} $(\mcal{X}(P^{\centerdot}), \mcal{Y}(P^{\centerdot}))$ is a torsion 
theory for $\cat{Mod}A$ and $\opn{H}^{-1}(P^{\centerdot}) \in 
\mcal{Y}(P^{\centerdot})$.
\end{enumerate}
\end{cor}

\begin{proof}  By Lemmas \ref{lem2.6}, \ref{lem2.7}
and Theorem \ref{thm3.5}.
\end{proof}

For an object $X$ in an additive category $\mcal{B}$,
we denote by $\cat{add}(X)$ the full subcategory of $\mcal{B}$
consisting of objects which are direct summands of finite
direct sums of copies of $X$.

\begin{cor} \label{cor3.7} 
For any tilting complexes 
$P^{\centerdot}_{1} : P^{-1}_{1} \arr P^{0}_{1}$, 
$P^{\centerdot}_{2} : P^{-1}_{2} \arr P^{0}_{2}$ for $A$
of term length two, the following are equivalent.
\begin{enumerate}
\rmitem{1} $(\mcal{X}(P^{\centerdot}_{1}), \mcal{Y}(P^{\centerdot}_{1}))
= (\mcal{X}(P^{\centerdot}_{2}), \mcal{Y}(P^{\centerdot}_{2}))$.
\rmitem{2} $\cat{add}(P^{\centerdot}_{1}) = \cat{add}(P^{\centerdot}_{2})$
in $\cat{K}^{\mrm{b}}(\cat{Proj} A)$.
\end{enumerate}
\end{cor}

\begin{proof}  
(1) $\Rightarrow$ (2).
It follows by Corollary \ref{cor3.6} that 
$Q^{\centerdot} = P^{\centerdot}_{1}\oplus P^{\centerdot}_{2}$ is a tilting
complex such that $(\mcal{X}(Q^{\centerdot}), \mcal{Y}(Q^{\centerdot}))
= (\mcal{X}(P^{\centerdot}_{i}), \mcal{Y}(P^{\centerdot}_{i}))$ ($i = 1,2$).
Let $B = \opn{End}_{\cat{D}(\cat{Mod}A)}(Q^{\centerdot})^{\mrm{op}}$ and for $i = 1,2$
denote by $e_{i}$ the composite of canonical homomorphisms $Q^{\centerdot} \arr 
P^{\centerdot}_{i} \arr Q^{\centerdot}$ .  Then for $i = 1,2$ we have an 
equivalence $\cat{D}^{-}(\cat{Mod}B) \arr \cat{D}^{-}(\cat{Mod}e_{i}Be_{i})$ 
which sends $Be_{i}$ to $e_{i}Be_{i}$, so that the $Be_{i}$ are tilting complexes 
for $B$, i.e. projective generators for $\cat{Mod}B$.  It follows by Morita 
Theory that $\cat{add}B = \cat{add}Be_{i}$ in $\cat{Mod}B$. Thus 
$\cat{add}(P^{\centerdot}_{1}) = 
\cat{add}(P^{\centerdot}_{2})$ in $\cat{K}^{\mrm{b}}(\cat{Proj} A)$.

(2) $\Rightarrow$ (1).
It is obviously deduced that
$\cat{add}(\opn{H}^{-1}(\nu(P^{\centerdot}_{1}))) =
\cat{add}(\opn{H}^{-1}(\nu(P^{\centerdot}_{2})))$ and
$\cat{add}(\opn{H}^{0}(P^{\centerdot}_{1}))
= \cat{add}(\opn{H}^{0}(P^{\centerdot}_{2}))$.
\end{proof}

\begin{thm} \label{thm3.8}
Assume $\mcal{X}(P^{\centerdot}) \cap
\mcal{Y}(P^{\centerdot}) = \{ 0 \}$ and $\opn{H}^{0}(P^{\centerdot}) \in 
\mcal{X}(P^{\centerdot})$.  Then the following hold.
\begin{enumerate}
\rmitem{1} $\{P^{\centerdot}[i] : i \in \mbb{Z}\}$ is a generating set
for $\cat{D}(\cat{Mod}A)$.
\rmitem{2} $\mcal{C}(P^{\centerdot})$ is admissible abelian.
\rmitem{3} $(\mcal{Y}(P^{\centerdot})[1], \mcal{X}(P^{\centerdot}))$ is a torsion 
theory for $\mcal{C}(P^{\centerdot})$.
\rmitem{4} The functor
\[
\opn{Hom}_{\cat{D}(\cat{Mod}A)}(P^{\centerdot}, -) : 
\mcal{C}(P^{\centerdot}) \arr \cat{Mod}B
\]
is an  equivalence.
\end{enumerate}
\end{thm}

\begin{proof}
By Lemma \ref{lem2.6} and Theorem \ref{thm2.15}.
\end{proof}

\begin{rem} \label{rem3.9}
The following are equivalent.
\begin{enumerate}
\rmitem{1} $P^{\centerdot}$ is a tilting complex.
\rmitem{2} $\mcal{X}(P^{\centerdot}) \cap \mcal{Y}(P^{\centerdot}) = \{ 0\}$
and $P^{\centerdot} \in \mcal{C}(P^{\centerdot})$.
\end{enumerate}
\end{rem}

\begin{exmp}[cf. \cite{HK}] \label{exa3.10}
Let $A$ be a finite dimensional algebra over a field  $k$
given by a quiver
\[ \begin{CD}
1 @>\alpha>>  2 \\
@A\delta AA @VV\beta V \\
4 @<<\gamma < 3
\end{CD}\]
with relations $\beta \alpha = \gamma \beta = \delta \gamma = \alpha \delta = 0$.
For each vertex $i$, we denote by $S(i), P(i)$ the corresponding simple and 
indecomposable projective left $A$-modules, respectively.  Define a complex 
$P^{\centerdot}$ as the mapping cone of the homomorphism
\[
d^{-1}_{P} =
\left[ \begin{smallmatrix}
f & 0 & 0 & 0 \\
0 & 0 & g & 0
\end{smallmatrix}\right]:
P(2)^2\oplus P(4)^2 
\arr
P(1)\oplus P(3) ,
\]
where $f$ and $g$ denote the right multiplications of $\alpha$ and $\gamma$, respectively.
Then $P^{\centerdot}$ is not a tilting complex.  However, $P^{\centerdot}$ satisfies the 
assumption of Theorem \ref{thm3.8} and hence we have an equivalence of abelian categories
\[
\opn{Hom}_{\cat{D}(\cat{Mod}A)}(P^{\centerdot}, -) : 
\mcal{C}(P^{\centerdot}) \arr \cat{Mod}B ,
\]
where $B= \opn{End}_{\cat{D}(\cat{Mod}A)}(P^{\centerdot})^{\mrm{op}}$ is a finite dimensional 
$k$-algebra given by a quiver
\[
1 \larr 2 \qquad  3 \larr 4 .
\]
There exist exact sequences in $\mcal{C}(P^{\centerdot})$ of the form
\[
0 \arr S(1) \arr S(2)[1] \arr P(1)[1] \arr 0,  \quad
0 \arr S(3) \arr S(4)[1] \arr P(3)[1] \arr 0,
\]
and these objects and morphisms generate $\mcal{C}(P^{\centerdot})$. 
\end{exmp}

\section{Equivalences between Torsion Theories} \label{s4}

Throughout this section, $P^{\centerdot}$ is assumed to be a tilting complex.  
Then there exists an equivalence of triangulated categories
\[     
F : \cat{D}^{-}(\cat{Mod}B) \arr \cat{D}^{-}(\cat{Mod}A)
\]
such that $F(B) = P^{\centerdot}$.  Let $G : \cat{D}^{-}(\cat{Mod}A) \arr
\cat{D}^{-}(\cat{Mod}B)$ be a quasi-inverse of $F$.
For any $n \in \mbb{Z}$, we have ring homomorphisms
\[
B \arr \opn{End}_{A}(\opn{H}^{n}(P^{\centerdot}))^{\mrm{op}} \quad\text{and}\quad
B \arr \opn{End}_{A}(\opn{H}^{n}((P^{\centerdot})^{*})).
\]
In particular, $\opn{H}^{0}(P^{\centerdot})$ is an $A$-$B$-bimodule and
$\opn{H}^{1}((P^{\centerdot})^{*})$ is a $B$-$A$-bimodule.

\begin{lem} \label{lem4.1}
The following hold.
\begin{enumerate}
\rmitem{1} For any $X^{\centerdot} \in \mcal{C}(P^{\centerdot})$, we have 
$G(X^{\centerdot}) \cong \opn{Hom}_{\cat{D}(\cat{Mod}A)}(P^{\centerdot}, X^{\centerdot})$.
\rmitem{2} We have an equivalence
\[
\opn{Hom}_{\cat{D}(\cat{Mod}A)}(P^{\centerdot}, -) : \mcal{C}(P^{\centerdot}) 
\arr \cat{Mod}B
\]
whose quasi-inverse is given by the restriction of $F$ to $\cat{Mod}B$.
\end{enumerate}
\end{lem}

\begin{proof}
See \cite{Ri}, Section 4.
\end{proof}

\begin{lem} \label{lem4.2}
There exists a tilting complex $Q^{\centerdot} \in \cat{K}^{\mrm{b}}(\cat{proj}B)$
such that
\begin{enumerate}
\rmitem{i} $Q^{\centerdot} \cong G(A)$,
\rmitem{ii} $Q^i = 0$ for $i > 1$ and $i < 0$,
\rmitem{iii} $\opn{H}^i(Q^{\centerdot}) \cong \opn{H}^i((P^{\centerdot})^{*})$
for $0 \leq i \leq 1$,
\rmitem{iv} $\opn{H}^i(\opn{Hom}^{\centerdot}_{B}(Q^{\centerdot}, B)) 
\cong \opn{H}^i(P^{\centerdot})$
for $- 1 \leq i \leq 0$.
\end{enumerate}
\end{lem}

\begin{proof}
By \cite{Ri}, Proposition 6.3, there exists $Q^{\centerdot} \in
\cat{K}^{\mrm{b}}(\cat{proj}B)$ satisfying $Q^{\centerdot} \cong G(A)$.
Since 
\[ \begin{aligned}
\opn{H}^i(Q^{\centerdot}) & \cong
\opn{Hom}_{\cat{D}(\cat{Mod}B)}(B, Q^{\centerdot}[i]) \\
& \cong
\opn{Hom}_{\cat{D}(\cat{Mod}A)}(P^{\centerdot}, A[i]) \\
& \cong \opn{H}^i((P^{\centerdot})^{*}),
\end{aligned} \]
we have $Q^{\centerdot} \cong {\sigma}_{\leq 1}(Q^{\centerdot})$ 
in $\cat{K}^{\mrm{b}}(\cat{proj}B)$.
Also, since \[ \begin{aligned}
\opn{H}^i(\opn{Hom}^{\centerdot}_{B}(Q^{\centerdot}, B)) & \cong 
\opn{Hom}_{\cat{D}(\cat{Mod}B)}(Q^{\centerdot}, B[i]) \\
& \cong \opn{Hom}_{\cat{D}(\cat{Mod}A)}(A, P^{\centerdot}[i]) \\
& \cong \opn{H}^i(P^{\centerdot}),
\end{aligned} \]
we have $\opn{Hom}^{\centerdot}_{B}(Q^{\centerdot}, B)
\cong {\sigma}_{\leq 0}(\opn{Hom}^{\centerdot}_{B}(Q^{\centerdot}, B))$
in $\cat{K}^{\mrm{b}}(\cat{proj}B^{\mrm{op}})$ and
$Q^{\centerdot} \cong {\sigma '}_{\geq 0}(Q^{\centerdot})$ 
in $\cat{K}^{\mrm{b}}(\cat{proj}B)$.
Thus, we can assume $Q^{i} = 0$ for $i > 1$ and $i < 0$.
\end{proof}

\begin{lem} \label{lem4.3}
For any $M \in \cat{Mod}B$, the following hold.
\begin{enumerate}
\rmitem{1} $\opn{H}^i(F(M)) = 0$ for $i > 0$ and $i < - 1$.
\rmitem{2} $\opn{H}^0(F(M)) \cong \opn{H}^0(P^{\centerdot}){\otimes}_{B}M$.
\rmitem{3} $\opn{H}^{-1}(F(M)) \cong 
\opn{Hom}_{B}(\opn{H}^{1}((P^{\centerdot})^{*}), M)$.
\end{enumerate} 
\end{lem}

\begin{proof}
For any $i \in \mbb{Z}$, we have
\[ \begin{aligned}
\opn{H}^i(F(M)) & \cong \opn{Hom}_{\cat{D}(\cat{Mod}A)}(A, F(M)[i]) \\
& \cong \opn{Hom}_{\cat{D}(\cat{Mod}B)}(Q^{\centerdot}, M[i]).
\end{aligned} \]
Thus $\opn{H}^i(F(M)) = 0$ for $i > 0$ and $i < - 1$.  Also,
\[ \begin{aligned}
\opn{H}^0(F(M)) & \cong \opn{Hom}_{\cat{D}(\cat{Mod}B)}(Q^{\centerdot}, 
M)  \\
& \cong \opn{H}^{0}(\opn{Hom}^{\centerdot}_{B}(Q^{\centerdot}, M)) \\
& \cong \opn{H}^{0}(\opn{Hom}^{\centerdot}_{B}(Q^{\centerdot}, B){\otimes}_{B}M) \\
& \cong \opn{H}^{0}(\opn{Hom}^{\centerdot}_{B}(Q^{\centerdot}, B)){\otimes}_{B}M \\
& \cong \opn{H}^{0}(P^{\centerdot}){\otimes}_{B}M,
\end{aligned}
\] \[
\begin{aligned}
\opn{H}^{- 1}(F(M)) & \cong
\opn{Hom}_{\cat{D}(\cat{Mod}B)}(Q^{\centerdot}, M[- 1]) \\
& \cong \opn{H}^{- 1}(\opn{Hom}^{\centerdot}_{B}(Q^{\centerdot}, M)) \\
& \cong \opn{Hom}_{B}(\opn{H}^{1}(Q^{\centerdot}), M) \\
& \cong \opn{Hom}_{B}(H^1((P^{\centerdot})^{*}), M).
\end{aligned} \]
\end{proof}

\begin{thm} \label{thm4.4}
Define a pair of full subcategories of $\cat{Mod}B$
\[
\mcal{U}(P^{\centerdot}) = \opn{Ker}(\opn{H}^0(P^{\centerdot}){\otimes}_{B}-), \quad
\mcal{V}(P^{\centerdot}) = \opn{Ker}(\opn{Hom}_{B}(\opn{H}^1((P^{\centerdot})^{*}), -)).
\]
Then the following hold.
\begin{enumerate}
\rmitem{1} $(\mcal{U}(P^{\centerdot}), \mcal{V}(P^{\centerdot}))$ 
is a torsion theory for $\cat{Mod}B$.
\rmitem{2} We have a pair of functors
\[
\opn{Hom}_{A}(\opn{H}^0(P^{\centerdot}), -) : \mcal{X}(P^{\centerdot})
\arr \mcal{V}(P^{\centerdot}), \quad
\opn{H}^0(P^{\centerdot}){\otimes}_{B}- : \mcal{V}(P^{\centerdot}) 
\arr \mcal{X}(P^{\centerdot})
\]
which define an equivalence.
\rmitem{3} We have a pair of functors
{\small \[
\opn{H}^1((P^{\centerdot})^{*}){\otimes}_{A}- : \mcal{Y}(P^{\centerdot}) 
\arr \mcal{U}(P^{\centerdot}), \quad
\opn{Hom}_{B}(\opn{H}^1((P^{\centerdot})^{*}), -) : \mcal{U}(P^{\centerdot}) 
\arr \mcal{Y}(P^{\centerdot})\]}
which define an equivalence.
\end{enumerate}
\end{thm}

\begin{proof}
(1) According to Lemmas \ref{lem3.2} and \ref{lem4.2}, we can apply Corollary \ref{cor3.6} for 
a tilting complex $Q^{\centerdot}$ to conclude that 
$(\mcal{U}(P^{\centerdot}), \mcal{V}(P^{\centerdot}))$ is a 
torsion theory for $\cat{Mod}B$.

(2) For any $X \in \mcal{X}(P^{\centerdot})$, by Lemmas \ref{lem2.13}, 
\ref{lem4.1} (1) and \ref{lem4.3} (3) we have
\[\begin{aligned}
\opn{Hom}_{B}(\opn{H}^{1}((P^{\centerdot})^{*}), 
\opn{Hom}_{A}(\opn{H}^{0}(P^{\centerdot}), X)) & \cong
\opn{H}^{-1}(F(G(X))) \\
& \cong \opn{H}^{-1}(X) \\
& = 0 .
\end{aligned}\]
Also, since by Lemma \ref{lem3.2} (1) and Corollary \ref{cor3.6} 
$\opn{H}^{1}((P^{\centerdot})^{*}){\otimes}_{A}
\opn{H}^{0}(P^{\centerdot}) = 0$, \\ $\opn{H}^{1}((P^{\centerdot})^{*}){\otimes}_{A}
\opn{H}^{0}(P^{\centerdot}){\otimes}_{B}M = 0$ for all $M \in \mcal{V}(P^{\centerdot})$.
The last assertion follows by Lemmas \ref{lem2.13}, \ref{lem4.1} and \ref{lem4.3}.

(3) For any $Y \in \mcal{Y}(P^{\centerdot})$, by Lemmas \ref{lem2.13}, \ref{lem3.1}, 
\ref{lem4.1} (1) and \ref{lem4.3} (2) we have
\[\begin{aligned}
\opn{H}^{0}(P^{\centerdot}){\otimes}_{B}
\opn{H}^{1}((P^{\centerdot})^{*}){\otimes}_{A}Y & \cong
\opn{H}^{0}(F(G(Y[1]))) \\
& \cong \opn{H}^{0}(Y[1]) \\
& = 0 .
\end{aligned}\]
Also, since $\opn{H}^{1}((P^{\centerdot})^{*}){\otimes}_{A}
\opn{H}^{0}(P^{\centerdot}) = 0$, for any $N \in \mcal{U}(P^{\centerdot})$ we have
\[\begin{aligned}
\opn{Hom}_{A}(\opn{H}^{0}(P^{\centerdot}), 
\opn{Hom}_{B}(\opn{H}^{1}((P^{\centerdot})^{*}), N))  & \cong
\opn{Hom}_{B}(\opn{H}^{1}((P^{\centerdot})^{*})
{\otimes}_{A}\opn{H}^{0}(P^{\centerdot}), N)  \\
& = 0.
\end{aligned}\]
The last assertion follows by Lemmas 
\ref{lem2.13}, \ref{lem4.1} and \ref{lem4.3}.
\end{proof}

\begin{defn} \label{dfn4.5}
Let $(\mcal{U}, \mcal{V})$ be a torsion theory for an abelian category 
$\mcal{A}$.  
Then $(\mcal{U}, \mcal{V})$ is called splitting if  
$\opn{Ext}^{1}_{\mcal{A}}(\mcal{V}, \mcal{U}) = 0$.
\end{defn}

For a left $A$-module $M$, we denote by $\opn{proj\ dim}{}_{A}M$
(resp., $\opn{inj\ dim}{}_{A}M$) the projective (resp., the injective) 
dimension of $M$.

\begin{prop} \label{prop4.6}
The torsion theory $(\mcal{U}(P^{\centerdot}), \mcal{V}(P^{\centerdot}))$ 
for $\cat{Mod}B$ is splitting
if and only if $\opn{Ext}^{2}_{A}(\mcal{X}(P^{\centerdot}),
\mcal{Y}(P^{\centerdot})) = 0$.
In particular, $(\mcal{U}(P^{\centerdot}), \mcal{V}(P^{\centerdot}))$ 
is splitting
if either $\opn{proj\ dim}X \leq 1$ for all $X \in \mcal{X}(P^{\centerdot})$ 
or $\opn{inj\ dim} Y \leq 1$ for all $Y \in \mcal{Y}(P^{\centerdot})$.
\end{prop}

\begin{proof}
For any $X \in \mcal{X}(P^{\centerdot})$ and $Y \in 
\mcal{Y}(P^{\centerdot})$, we have 
{\small \[ \begin{aligned}
\opn{Ext}^{1}_{B}(\opn{Hom}_{A}(\opn{H}^{0}(P^{\centerdot}),X), 
\opn{H}^{1}((P^{\centerdot})^{*}){\otimes}_{A}Y)
& \cong \opn{Hom}_{\cat{D}(\cat{Mod}B)}(G(X), G(Y[1])[1]) \\
& \cong \opn{Hom}_{\cat{D}(\cat{Mod}A)}(X,Y[2]) \\
& \cong \opn{Ext}^{2}_{A}(X,Y).  
\end{aligned} \]}
\end{proof}

\section{Torsion Theories for Artin Algebras} \label{s5}

In this section, we deal with the case where $R$ is a commutative artin ring, 
$I$ is an injective envelope of an $R$-module $R/\opn{rad}(R)$ and $A$ is a 
finitely generated $R$-module.
We denote by $\cat{mod}A$ the full abelian subcategory of 
$\cat{Mod}A$ consisting of finitely generated modules.
Note that $\opn{H}^{n}(P^{\centerdot}),  \opn{H}^{n}(\nu(P^{\centerdot})) \in \cat{mod}A$
for all $n \in \mbb{Z}$.  We set
\[     
\mcal{X}_{c}(P^{\centerdot}) = \mcal{X}(P^{\centerdot}) \cap \cat{mod}A     
\quad\text{and}\quad
\mcal{Y}_{c}(P^{\centerdot}) = \mcal{Y}(P^{\centerdot}) \cap \cat{mod}A.
\]

\begin{prop} \label{prop5.1}
The following are equivalent.
\begin{enumerate}
\rmitem{1} $\mcal{X}_{c}(P^{\centerdot}) \cap \mcal{Y}_{c}(P^{\centerdot}) = \{ 0\}$
and $\opn{H}^{0}(P^{\centerdot}) \in \mcal{X}_{c}(P^{\centerdot})$.
\rmitem{2} $\mcal{X}_{c}(P^{\centerdot}) \cap \mcal{Y}_{c}(P^{\centerdot}) = \{ 0\}$
and $\tau(X) \in \mcal{X}_{c}(P^{\centerdot})$, $\pi(X) \in 
\mcal{Y}_{c}(P^{\centerdot})$ for all $X \in \cat{mod}A$.
\rmitem{3} $(\mcal{X}_{c}(P^{\centerdot}), \mcal{Y}_{c}(P^{\centerdot}))$ is a torsion 
theory for $\cat{mod}A$.
\rmitem{4} $\mcal{X}_{c}(P^{\centerdot})$ consists of the modules generated by 
$\opn{H}^{0}(P^{\centerdot})$ and $\mcal{Y}_{c}(P^{\centerdot})$ consists of the modules 
cogenerated by $\opn{H}^{-1}(\nu(P^{\centerdot}))$.
\end{enumerate}
\end{prop}
     
\begin{proof}          
By the same arguments as in the proof of Theorem \ref{thm3.5}.
\end{proof}

\begin{lem} \label{lem5.2}
The following are equivalent.
\begin{enumerate}
\rmitem{1} $\{P^{\centerdot}[i] : i \in \mbb{Z}\}$ is a generating set
for $\cat{D}(\cat{mod}A)$.
\rmitem{2} $\mcal{X}_{c}(P^{\centerdot}) \cap \mcal{Y}_{c}(P^{\centerdot}) = \{0 \}$.
\end{enumerate}
\end{lem}

\begin{proof}          
By the same arguments as in the proof of Lemma \ref{lem2.6}.
\end{proof}

\begin{lem} \label{lem5.3}
The following hold.
\begin{enumerate}
\rmitem{1} If $DA \in \mcal{X}_{c}(P^{\centerdot})$, then 
$\opn{H}^{-1}(P^{\centerdot}) = 0$, i.e. 
$P^{\centerdot} \cong \opn{H}^{0}(P^{\centerdot})$ in $\cat{D}(\cat{mod}A)$.
\rmitem{2} $\opn{H}^{0}(\nu(P^{\centerdot})) \in \mcal{X}_{c}(P^{\centerdot})$
if and only if $\opn{H}^{-1}(P^{\centerdot}) \in \mcal{Y}_{c}(P^{\centerdot})$.
\end{enumerate}
\end{lem}

\begin{proof}          
For any $P \in \cat{proj}A$, we have 
functorial isomorphisms
\[     
\nu(P) \cong DA{\otimes}_{A}P \quad\text{and}\quad P \cong
\opn{Hom}_{A}(DA, \nu(P)).
\]
Thus 
\[
\opn{H}^{0}(\nu(P^{\centerdot})) \cong DA{\otimes}_{A}
\opn{H}^{0}(P^{\centerdot}) \quad\text{and}\quad \opn{H}^{-1}(P^{\centerdot}) \cong 
\opn{Hom}_{A}(DA, \opn{H}^{-1}(\nu(P^{\centerdot})))
\] and hence
\[ \begin{aligned}
\opn{Hom}_{A}(\opn{H}^{0}(\nu(P^{\centerdot})), \opn{H}^{-1}(\nu(P^{\centerdot})))
& \cong 
\opn{Hom}_{A}(DA {\otimes}_{A} \opn{H}^{0}(P^{\centerdot}),
\opn{H}^{-1}(\nu(P^{\centerdot}))) \\
& \cong \opn{Hom}_{A}(\opn{H}^{0}(P^{\centerdot}), \opn{Hom}_{A}(DA , 
\opn{H}^{-1}(\nu(P^{\centerdot})))) \\
& \cong \opn{Hom}_{A}(\opn{H}^{0}(P^{\centerdot}), \opn{H}^{-1}(P^{\centerdot})).
\end{aligned} \] \end{proof}

\begin{lem} \label{lem5.4}
Assume $\mcal{X}_{c}(P^{\centerdot}) \cap \mcal{Y}_{c}(P^{\centerdot}) = \{ 0\}$ and 
$\opn{H}^{0}(P^{\centerdot}) \in \mcal{X}_{c}(P^{\centerdot})$.  
Then the following are equivalent.
\begin{enumerate}
\rmitem{1} $\opn{H}^{0}(\nu(P^{\centerdot})) \in \mcal{X}_{c}(P^{\centerdot})$.
\rmitem{2} $\mcal{X}_{c}(P^{\centerdot})$ is stable under $DA
{\otimes}_{A} -$.
\rmitem{3} $\opn{H}^{-1}(P^{\centerdot}) \in \mcal{Y}_{c}(P^{\centerdot})$.
\rmitem{4} $\mcal{Y}_{c}(P^{\centerdot})$ is stable under $\opn{Hom}_{A}(DA, -)$.
\end{enumerate}
\end{lem}

\begin{proof}
(1) $\Rightarrow$ (2).  Let $X \in \mcal{X}_{c}(P^{\centerdot})$.
Then by Proposition \ref{prop5.1} $X$ is generated by 
$\opn{H}^{0}(P^{\centerdot})$ and hence $DA{\otimes}_{A}X$ is generated 
by $DA {\otimes}_{A} \opn{H}^{0}(P^{\centerdot}) \cong 
\opn{H}^{0}(\nu(P^{\centerdot})) \in \mcal{X}_{c}(P^{\centerdot})$.

(2) $\Rightarrow$ (3).  Since $\opn{H}^{0}(\nu(P^{\centerdot})) \cong 
DA {\otimes}_{A} \opn{H}^{0}(P^{\centerdot}) \in 
\mcal{X}_{c}(P^{\centerdot})$, by Lemma \ref{lem5.3} (2) we have 
$\opn{H}^{-1}(P^{\centerdot}) \in \mcal{Y}_{c}(P^{\centerdot})$.

(3) $\Rightarrow$ (4) $\Rightarrow$ (1).  By the dual arguments.
\end{proof}

\begin{prop} \label{prop5.5}
The following are equivalent.
\begin{enumerate}
\rmitem{1} $P^{\centerdot}$ is a tilting complex.
\rmitem{2} $\mcal{X}_{c}(P^{\centerdot}) \cap \mcal{Y}_{c}(P^{\centerdot}) = \{ 0\}$,
$\opn{H}^{0}(P^{\centerdot}) \in \mcal{X}_{c}(P^{\centerdot})$
and $\opn{H}^{-1}(P^{\centerdot}) \in \mcal{Y}_{c}(P^{\centerdot})$.
\rmitem{3} $(\mcal{X}_{c}(P^{\centerdot}), \mcal{Y}_{c}(P^{\centerdot}))$ 
is a torsion theory for $\cat{mod}A$ and 
$\opn{H}^{-1}(P^{\centerdot})$ $\in \mcal{Y}_{c}(P^{\centerdot})$.
\rmitem{4} $(\mcal{X}_{c}(P^{\centerdot}), \mcal{Y}_{c}(P^{\centerdot}))$
is a torsion theory for $\cat{mod}A$ and 
$\mcal{X}_{c}(P^{\centerdot})$ is stable under $DA {\otimes}_{A} -$.
\rmitem{5} $(\mcal{X}_{c}(P^{\centerdot}), \mcal{Y}_{c}(P^{\centerdot}))$ is 
a torsion theory for $\cat{mod}A$ and 
$\mcal{Y}_{c}(P^{\centerdot})$ is stable under $\opn{Hom}_{A}(DA, -)$.
\end{enumerate}
\end{prop}

\begin{proof}
By Proposition \ref{prop5.1}, Lemmas \ref{lem2.7},
\ref{lem5.2} and \ref{lem5.4}.
\end{proof}

\begin{defn} \label{dfn5.6}
Let $\mcal{A}$ be an abelian category and $\mcal{C}$ a full 
subcategory of $\mcal{A}$ closed under extensions.  
Then an object $X \in \mcal{C}$ is called $\opn{Ext}$-projective 
(resp., $\opn{Ext}$-injective) if
$\opn{Ext}^{1}_{\mcal{A}} (X, \mcal{C}) = 0$ 
(resp., $\opn{Ext}^{1}_{\mcal{A}} (\mcal{C}, X) = 0$). 
\end{defn}

\begin{prop} \label{prop5.7}
Assume $P^{\centerdot}$ is a tilting complex.  Then the following hold.
\begin{enumerate}
\rmitem{1} $\opn{H}^{0}(P^{\centerdot}) \in \mcal{X}_{c}(P^{\centerdot})$ is 
$\opn{Ext}$-projective and generates $\mcal{X}_{c}(P^{\centerdot})$.
\rmitem{2} $\opn{H}^{-1}(\nu(P^{\centerdot})) \in \mcal{Y}_{c}(P^{\centerdot})$ is 
$\opn{Ext}$-injective and cogenerates $\mcal{Y}_{c}(P^{\centerdot})$.
\end{enumerate}
\end{prop}

\begin{proof}
By Propositions \ref{prop5.1}, \ref{prop5.5} and 
Lemmas \ref{lem3.3}, \ref{lem3.4}.
\end{proof}

\begin{thm} \label{thm5.8}
Let $(\mcal{X} , \mcal{Y})$ be a torsion theory for $\cat{mod}A$ such 
that $\mcal{X}$ contains an $\opn{Ext}$-projective module X which generates 
$\mcal{X}$, $\mcal{Y}$ contains an $\opn{Ext}$-injective module Y which 
cogenerates $\mcal{Y}$, and $\mcal{X}$ is stable under $DA {\otimes}_{A} -$.
Let $M^{\centerdot}_{X}$ be a minimal projective presentation of $X$ and
$N^{\centerdot}_{Y}$ a minimal injective presentation of $Y$.
Then 
\[     
P^{\centerdot} = M^{\centerdot}_{X} \oplus 
\opn{Hom}^{\centerdot}_{A}(DA, N^{\centerdot}_{Y})[1]
\]
is a tilting complex such that 
$\mcal{X} = \mcal{X}_{c}(P^{\centerdot})$ and $\mcal{Y} = 
\mcal{Y}_{c}(P^{\centerdot})$.
\end{thm}

\begin{proof}
According to Proposition \ref{prop5.5}, we have only to show that 
$\mcal{X} = \mcal{X}_{c}(P^{\centerdot})$ and $\mcal{Y} = 
\mcal{Y}_{c}(P^{\centerdot})$.  
It follows by \cite{Ho2}, Lemmas 2 and 3 that
$\opn{H}^{0}(P^{\centerdot}) \in \mcal{X}$ and
$\opn{H}^{-1}(\nu(P^{\centerdot})) \in \mcal{Y}$.
Since $X$ is a direct summand of 
$\opn{H}^{0}(P^{\centerdot})$ and $Y$ is a direct summand of
$\opn{H}^{-1}(\nu(P^{\centerdot}))$, it follows that 
$\opn{H}^{0}(P^{\centerdot})$ generates $\mcal{X}$ and
$\opn{H}^{-1}(\nu(P^{\centerdot}))$ cogenerates $\mcal{Y}$.  
It now follows by Remark \ref{rem2.9}, Lemmas \ref{lem3.2}, \ref{lem3.3} (2) that $\mcal{X} = 
\mcal{X}_{c}(P^{\centerdot})$ and $\mcal{Y} = \mcal{Y}_{c}(P^{\centerdot})$.
\end{proof}


\end{document}